\pdfoutput=1
\RequirePackage{ifpdf}
\ifpdf % We~are running pdfTeX in pdf mode
\documentclass[pdftex]{sigma}
\else
\documentclass{sigma}
\fi

\usepackage{dsfont,esint}

\numberwithin{equation}{section}

\newtheorem{Theorem}{Theorem}[section]
\newtheorem{Corollary}[Theorem]{Corollary}
\newtheorem{Conjecture}[Theorem]{Conjecture}
\newtheorem{oproblem}[Theorem]{Open Problem}
 { \theoremstyle{definition}
\newtheorem{Definition}[Theorem]{Definition}
\newtheorem{Example}[Theorem]{Example}
\newtheorem{Remark}[Theorem]{Remark} }

\newcommand{\Rm}{{\rm Rm}}

\newcommand{\Vol}{{\rm Vol}}

\newcommand{\cH}{\mathcal{H}}

\newcommand{\cR}{\mathcal{R}}
\newcommand{\cS}{\mathcal{S}}

\begin{document}
\allowdisplaybreaks

\renewcommand{\thefootnote}{}

\renewcommand{\PaperNumber}{104}

\FirstPageHeading

\ShortArticleName{Conjectures and Open Questions on the Structure and Regularity of Spaces}

\ArticleName{Conjectures and Open Questions\\ on the Structure and Regularity of Spaces\\ with Lower Ricci Curvature Bounds\footnote{This paper is a~contribution to the Special Issue on Scalar and Ricci Curvature in honor of Misha Gromov on his 75th Birthday. The full collection is available at \href{https://www.emis.de/journals/SIGMA/Gromov.html}{https://www.emis.de/journals/SIGMA/Gromov.html}}}

\Author{Aaron NABER}

\AuthorNameForHeading{A.~Naber}

\Address{Department of Mathematics, Northwestern University, USA}
\Email{\href{mailto:anaber@math.northwestern.edu}{anaber@math.northwestern.edu}}

\ArticleDates{Received July 29, 2020, in final form October 11, 2020; Published online October 20, 2020}

\Abstract{In this short note we review some known results on the structure and regularity of spaces with lower Ricci curvature bounds. We present some known and new open questions about next steps.}

\Keywords{Ricci curvature; regularity}

\Classification{53C21; 53C23}

\section{Introduction to Ricci curvature and limits}\label{s:intro}

\subsection{What is Ricci curvature?}

The Riemannian curvature tensor $\Rm$ of a Riemannian manifold $\big(M^n,g\big)$ is a four tensor, that is it takes three vector fields and gives a vector field. Classically it is defined as the antisymmetric part of the Hessian:
\begin{gather*}
\Rm(X,Y)Z \equiv \nabla^2_{X,Y}Z - \nabla^2_{Y,X}Z ,	
\end{gather*}
however from an analytic point of view one should interpret this four tensor $\Rm$ as the Hessian of the Riemannian metric~$g$. The actual Hessian of $g$ is zero, essentially by definition, however $\Rm$ is the correct moral replacement. In particular, one should expect that a manifold with bounded curvature should behave much like a function on ${\mathbb R}^n$ with~$C^2$ bounds. Up to a choice of coordinate chart, this is indeed the case.

The Ricci curvature is then defined as the two tensor obtained by tracing
\begin{gather*}
{\rm Ric}(X,Y)\equiv \sum_{E_a}\langle\Rm(X,E_a)E_a,Y\rangle ,	
\end{gather*}
where $E_a$ is an orthonormal basis. If $\Rm$ is interpreted as the Hessian of the metric $g$, then one should interpret the Ricci curvature as the Laplacian of the Riemannian metric. Clearly, it is a highly nonlinear Laplacian and this moral only holds up to the diffeomorphism gauge group, but it still puts into perspective why the Ricci curvature should play a central role in so many situations. One can then similarly define from this the scalar curvature as the trace $R\equiv \sum_{E_a} {\rm Ric}(E_a,E_a)$.

\subsection{Limiting under lower Ricci curvature}

The study of the structure and regularity of spaces with lower bounds on Ricci curvature is an old topic. The possibility for making a systematic study of the structure of such spaces began with the proof of a compactness theorem:

\begin{Theorem}[Gromov compactness theorem] Let $\big(M^n_i,g_i,p_i\big)$ be a sequence of pointed Riemannian manifolds such that ${\rm Ric}_i\geq -\lambda$ uniformly. Then there exists a metric space $(X,d,p)$ with $p\in X$ such that, after possibly passing to a subsequence, we have
	\begin{gather*}
	\big(M^n,g_i,p_i\big)\to (X,d,p) ,	
	\end{gather*}
where the convergence is in the Gromov--Hausdorff topology.
\end{Theorem}

This compactness theorem immediately begs the question of what type of metric spaces~$X$ might exist as limits. Gromov's theorem gives only that $X$ is a metric space, without giving any further refined structure on the space\footnote{In fact, Gromov did also prove that $X$ is a length space.}. However, one should expect much more in principle. Though at this point one can really refine the following categories much more, we can roughly break our interest in limits
\begin{gather*}
	\big(M^n,g_i,p_i\big)\to (X,d,p) ,
\end{gather*}
into two general groups:
\begin{itemize}\itemsep=0pt
\item (lower Ricci noncollapsed) $\big(M^n_i,g_i,p_i\big)$ satisfy ${\rm Ric}_i\geq -\lambda$ and noncollapsing $\Vol(B_1(p_i))>v>0$,	
\item (lower Ricci collapsed) $\big(M^n_i,g_i,p_i\big)$ satisfy ${\rm Ric}_i\geq -\lambda$ and collapsing $\Vol(B_1(p_i))\to 0$.	
\end{itemize}

In each of these cases, one want to explore both the structure of the metric spaces $(X,d)$ and a~priori regularity of the sequence $\big(M^n_i,g_i\big)$. We will cover what is known in these cases and present some open questions on future directions.

\section{Lower Ricci curvature and noncollapsing}

Let us discuss in this section the structure and regularity of limits under a uniform lower Ricci curvature and noncollapsing assumption:
\begin{gather*}%\label{e:lowerRicciNoncollapsing}
\big(M^n_i,g_i,p_i\big)\to (X,d,p) \qquad \text{s.t.} \quad {\rm Ric}_i\geq -\lambda \quad \text{and} \quad \Vol(B_1(p_i))>v>0 .	
\end{gather*}

The first real structure theorem for even the regular part of such spaces was presented by Cheeger and Colding:

\begin{Theorem}[Cheeger--Colding \cite{ChC2}]\label{t:lrnc_manifold}
	Let $\big(M^n_i,g_i,p_i\big)\to (X,d,p)$ satisfy
\[ {\rm Ric}_i\geq -\lambda \qquad \text{and} \qquad \Vol(B_1(p_i))>v>0,
\] then $X$ is bi-H\"older to a manifold away from a set of codimension two.
\end{Theorem}

The proof of the above is based on a Federer type stratification theory, which we review in Section~\ref{ss:stratification}, as well as a Reifenberg type theorem in order to deal with the regularity of the regular part. The Reifenberg argument necessarily only gives bi-H\"older control over the regular set, and the first open question we mention is on whether this can be improved:

\begin{oproblem}Let $\big(M^n_i,g_i,p_i\big)\to (X,d,p)$ satisfy ${\rm Ric}_i\geq -\lambda$ and $\Vol(B_1(p_i))>v>0$, then is $X$ bi-Lipschitz to a manifold away from a set of codimension two?
\end{oproblem}

This has turned out to be a subtle question. In order to get some sense of this let us recall the following example:

\begin{Example}[Colding--Naber \cite{CoNa3}]
There exists a limit space $X$, which is homeomorphic to ${\mathbb R}^n$ and for which every tangent cone, see Section~\ref{ss:tangent_cones}, is ${\mathbb R}^n$. There is a distinguished point $p\in X$ such that for any $x\in X$ if $\gamma_x\colon [0,1]\to X$ is the geodesic from $p$, then for every pair $x,y\in X$ and any $\theta\in [0,2\pi)$ we can find $t_i\to 0$ such that the angle between $\gamma_x(t_i)$ and $\gamma_y(t_i)$ converges to $\theta$. As a consequence, if $u\colon B_1(p)\to {\mathbb R}^n$ is a homeomorphism onto its image with $u$ either harmonic or a collection of distance functions, then $u$ cannot be bi-Lipschitz.
\end{Example}

The above example tells us that the maps one would typically attempt to use to build this homeomorphism cannot be bi-Lipschitz.

The other remaining question is on the size of the set on which $X$ is a manifold:

\begin{Conjecture}\label{conj:codim3_manifold}	Let $\big(M^n_i,g_i,p_i\big)\to (X,d,p)$ satisfy ${\rm Ric}_i\geq -\lambda$ and $\Vol(B_1(p_i))>v>0$, then~$X$ is homeomorphic to a manifold away from a set of codimension three.
\end{Conjecture}

The challenge in the above conjecture is that if one only removes a set of codimension three, then you must contend with some non-small singular points. These singularities should still all be of the form of a cone over a sphere on this set, and hence homeomorphic to a manifold, but this is still an open and challenging question.

\subsection{Tangent cones of noncollapsed spaces}\label{ss:tangent_cones}

In order to discuss more refined structure of noncollapsed limits $X$ we need to build the stratified singular set. The two main ingredients in this is the introduction of tangent cones and of symmetries of these tangent cones. Let us begin with a definition:

\begin{Definition}[tangent cones]Given a metric space $(X,d)$ with $x\in X$, we say another metric space $X_x$ is a tangent cone of $X$ at $x\in X$	if there exists $r_j\to 0$ such that $\big(X,r_{j}^{-1}d,x\big)\to X_x$.
\end{Definition}

Thus a tangent cone is representative of the infinitesimal behavior of a metric space $X$ near the point $x$. The tangent cone $X_x$ is obtained by zooming onto smaller and smaller balls around~$x$ and taking a~limit.

It follows from Gromov's compactness theorem that tangent cones exist for Ricci limits, a result which holds whether or not the limit is noncollapsed. One of the most important contributions of the works of Cheeger--Colding was the proof of the following:

\begin{Theorem}[tangent cones are metric cones \cite{ChC1}]\label{t:tangentcone_metriccone}
Let $\big(M^n_i,g_i,p_i\big)\to (X,d,p)$ satisfy ${\rm Ric}_i\geq -\lambda$ and $\Vol(B_1(p_i))>v>0$, then every tangent cone $X_x$ is a~metric cone $X_x\equiv C(Y_x)$ over $Y_x$ with $\operatorname{diam} Y_x\leq \pi$.
\end{Theorem}
The study of a stratification theory for essentially all nonlinear equations\footnote{Nonlinear harmonic maps, minimal surfaces, Yang--Mills, etc.} begins with a~statement analogous to the above. In most contexts the proof of the above statement is actually quite simple, however it takes a good amount of deep analysis in the context of lower Ricci bounds.

It is known that tangent cones need not be unique at a given point. It follows from the stratification theory that away from a set of codimension two {\it every} tangent cone is~${\mathbb R}^n$, and hence unique. A conjecture from~\cite{CoNa2} is that this uniqueness should still hold down to a set of codimension three:

\begin{Conjecture}\label{conj:codim3_unique}Let $\big(M^n_i,g_i,p_i\big)\to (X,d,p)$ satisfy ${\rm Ric}\geq -(n-1)$ and $\Vol(B_1(p_i))>v>0$, then the tangent cones of $X$ are unique away from a set of codimension three.
\end{Conjecture}

To prove the above means controlling the top stratum of the singular set to a~much larger degree than currently exists. Recently it was proved in~\cite{ChJiNa} that tangent cones are at least unique away from a~set of $n-2$ measure zero.\footnote{This was more a consequence of a broader structure theory on the rectifiability of the singular set.} See Conjecture~\ref{conj:Sk_all_symmetries} on further refinements of the above conjecture.

Besides tangent cones being nonunique, examples in \cite{CoNa2} even provided examples where tangent cones are not homeomorphic:

\begin{Example}[Colding--Naber]	There exists a limit space~$X^5$ with $p\in X$ such that tangent cones at $p$ are homeomorphic to both ${\mathbb R}^5$ and $C\big({\mathbb C} P^2\#\overline{{\mathbb C} P^2} \big)$, depending on the blow up sequence.
\end{Example}

It was conjectured this should not happen often:

\begin{Conjecture}Let $\big(M^n_i,g_i,p_i\big)\to (X,d,p)$ satisfy ${\rm Ric}\geq -(n-1)$ and $\Vol(B_1(p_i))>v>0$, then the tangent cones of $X$ are all homoemorphic at a fixed point away from a set of codimension five.
\end{Conjecture}

In order to build the stratification of the singular set, we need to also discuss the symmetries of tangent cones:

\begin{Definition}[symmetries]We say a metric space $X$ is $0$-symmetric if it is a cone space $X\equiv C(Y)$. We say $X$ is $k$-symmetric if $X\equiv {\mathbb R}^k\times C(Y)$ additionally splits off	$k$-Euclidean factors.
\end{Definition}

\subsection{Stratification of the singular set}\label{ss:stratification}

Now we can introduce the stratification of a noncollapsed Ricci limit space $X$, in the sense of Federer:

\begin{Definition}[stratification]
We define the $k$-stratum $\cS^k(X)\subseteq X$ by
\begin{gather*}
\cS^k(X)\equiv\big\{x\in X\colon \text{no tangent cone at $x$ is $k+1$-symmetric}\big\} .	
\end{gather*}
\end{Definition}

Note the subtlety in the above definition, which goes back to Federer. If $x\not\in \cS^k(X)$, then this guarantees the existence of at least one tangent cone with $k+1$-degrees of symmetry. We have that this is a nested sequence
\begin{gather*}
\cS^0(X)\subseteq \cS^1(X)\subseteq\cdots\subseteq \cS^{n-1}(X)\subseteq X .	
\end{gather*}
\looseness=-1 By applying Theorem~\ref{t:tangentcone_metriccone} and the dimension reduction technique of Federer, \cite{ChC2} was able to show
\begin{gather*}
\dim \cS^k\leq k .	
\end{gather*}

It is an important additional result of \cite{ChC2} that $\cS^{n-1}(X)=\cS^{n-2}(X)$, which is the starting point for Theorem~\ref{t:lrnc_manifold}. In terms of additional structure for the stratification of singular sets, the main result is

\begin{Theorem}[Cheeger--Jiang--Naber \cite{ChJiNa}]\label{t:Sk_rectifiable}
Let $\big(M^n_i,g_i,p_i\big)\to (X,d,p)$ satisfy ${\rm Ric}\geq -\lambda$ and $\Vol(B_1(p_i))>v>0$, then
\begin{enumerate}\itemsep=0pt
\item[$1)$] $\cS^k(X)$ is $k$-rectifiable,
\item[$2)$] For $\cH^k$-a.e.\ $x\in \cS^k(X)$ we have that every tangent cone is $k$-symmetric.
\end{enumerate}
\end{Theorem}
\begin{Remark}
In order to prove the above one has to work with the quantitative stratification~$\cS^k_\epsilon$, first introduced in~\cite{CheegerNaber_Ricci}, which decomposes the stratification into nicer pieces.
In particular, one can show that~$\cS^k_\epsilon$ have finite $k$-dimensional measure.	
\end{Remark}

Proving the above involves going beyond the dimension reduction technique of Federer. A~consequence of the second statement, applied to the top stratum, is that away from a set of $n-2$ measure zero 'every' tangent cone is of the form ${\mathbb R}^{n-2}\times C\big(S^1_r\big)$, where~$S^1_r$ is a circle of radius $r\leq 1$.

The first statement, that $\cS^k$ is $k$-rectifiable, says that $\cS^k(X)$ has a $k$-manifold structure away from a set of measure zero. This 'away from a set of measure zero' statement turns out to be sharp as one can build the following example:

\begin{Example}[Li--Na \cite{LiNa}]\label{ex:k_cantor}
For $k\leq n-2$ there exists a limit $\big(M^n_i,g_i,p_i\big)\to (X,d,p)$ where $\sec_i\geq 0$ and $\Vol(B_1(p_i))>v>0$ such that $\cS^k$ is a $k$-rectifiable, $k$-cantor set.
\end{Example}

The example tells us that an open manifold structure on $\cS^k$ is not possible for $k\leq n-2$. However, one might ask the following:

\begin{Conjecture}\label{conj:Sk_manifold}
Let $\big(M^n_i,g_i,p_i\big)\to (X,d,p)$ satisfy ${\rm Ric} \geq -\lambda$ and $\Vol(B_1(p_i))>v>0$, then~$\cS^k$ is contained in a union of $k$-dimensional submanifolds. More precisely, there exist a~countable collection of bi-H\"older maps $\varphi_i\colon B_1\big(0^k\big)\to M$ such that if $\hat \cS^k(X)\equiv \bigcup_i \varphi_i\big(B_1\big(0^k\big)\big)$ then
\begin{enumerate}\itemsep=0pt
\item[$1)$] $\cS^k(X)\subseteq \hat \cS^k(X)$,
\item[$2)$] for each $x\in \hat \cS^k(X)$ we have that every tangent cone is at least $k$-symmetric.	
\end{enumerate}
\end{Conjecture}

Note it is important that we only have an inclusion in $(1)$ above as opposed to equality by Example~\ref{ex:k_cantor}. A consequence of the above conjecture is the following:

\begin{Conjecture}\label{conj:Sk_all_symmetries}
Let $\big(M^n_i,g_i,p_i\big)\to (X,d,p)$ satisfy ${\rm Ric}\geq -\lambda$ and $\Vol(B_1(p_i))>v>0$, then for each $k\leq n-2$ there exists $\hat\cS^k$ with $\dim\hat\cS^k\leq k$ such that for each $x\in X\setminus \hat\cS^k$ we have that {\it every} tangent cone is at least $k+1$-symmetric.
\end{Conjecture}

Note that a positive answer to the above conjecture would solve both Conjecture \ref{conj:codim3_unique} and be a major step toward Conjecture~\ref{conj:codim3_manifold}. If the $M_i$ are K\"ahler then the above two conjectures have been solved by the very nice result in~\cite{LiSz}. Indeed, in this case one has equality $\cS^k(X)\equiv \hat \cS^k(X)$ as the complex analytic nature of the problem forces the singularities to be varieties.\footnote{a structure which is not possible in the real case by Example \ref{ex:k_cantor}.}

\subsection{Regularity of noncollapsed lower Ricci spaces and the energy identity}\label{ss:energy_identity}

Estimates and apriori estimates are often times a major goal in the study of solutions of equations. On noncollapsed spaces with lower Ricci bounds the best one has at the moment is the following:

\begin{Theorem}[Ji--Na \cite{JiNa}]	Let $\big(M^n,g,x\big)$ satisfy ${\rm Ric}\geq -\lambda$ and $\Vol(B_1(x))>v>0$, then for each $0<p<1$ there exists $C(n,\lambda,v,p)>0$ such that
\begin{gather*}
\fint_{B_1(x)} |{\rm Ric}|^p < C .	
\end{gather*}
\end{Theorem}

The above estimate seems less than optimal, and one might predict that an $L^1$ estimate on the Ricci curvature is possible, which is the context of Yau's conjecture:

\begin{Conjecture}[noncollapsing Yau's conjecture]
Let $\big(M^n,g,x\big)$ satisfy
\[ {\rm Ric}\geq -\lambda \qquad \text{and} \qquad \Vol(B_1(p))>v>0,\] then $\exists$ $C(n,\lambda,v)>0$ such that
\begin{gather*}
\fint_{B_1(x)} R < C .	
\end{gather*}
\end{Conjecture}

Recall that if there is a lower bound on the Ricci curvature, then an~$L^1$ bound on the Ricci and scalar curvature are equivalent. See Conjecture~\ref{conj:yau} for a more general version of the above.

Once an $L^1$ estimate on the scalar curvature is established, one can start to ask questions about the behavior of the scalar curvature $R\, {\rm d}v_g$ as a measure. The best way to understand this is to again consider sequences $\big(M^n_j,g_j,x_j\big)\to (X,d,x)$ and let the measures
\begin{gather*}
 R_j\,{\rm d}v_{g_j}\to \mu ,	
\end{gather*}
converge as measures to a limit on~$X$. To understand the behavior of $\mu$ let mention an example:

\begin{Example}Let $X=C\big(S^1_r\big)$ be a~two dimensional ice cream cone, which metrically is a cone over a~circle of radius $r<1$. Using warped coordinates we can construct smooth manifolds~$M_j$ which are isometric to~$X$ outside a ball $B_{j^{-1}}(x_j)$ around the cone point, and have globally nonnegative sectional curvature. It is clear that $R_j=0$ outside $B_{j^{-1}}(x_j)$, and from the Gauss--Bonnet formula we can compute that $\int_{M_j} R_j = 2\pi(1-r)$. In particular, $R_j \, {\rm d}v_{g_j}\to 2\pi(1-r) \delta_0$ converges to a Dirac delta measure at the cone point.
\end{Example}

From the above example we see that we should expect the scalar curvature measure to concentrate along the top stratum of the singular set. Conjecturally, this is all that should happen. To be precise recall from Theorem \ref{t:Sk_rectifiable} that away from a~set of $n-2$ measure zero there is a radius function $r_x$ such that at any~$x$ in the domain we have that the tangent cone at~$x$ is unique and isometric to ${\mathbb R}^{n-2}\times C\big(S^1_{r_x}\big)$. In reference to other areas, particularly Yang--Mills, we call the following the energy identity conjecture:

\begin{Conjecture}[energy identity]\label{con:2.20}
Let $(M^n_j,g_j,x_j)\to (X,d,x)$ satisfy ${\rm Ric}_j\geq -\lambda$ and the noncollapsing condition $\Vol(B_1(x_j))\geq v>0$, then as measures the scalar curvature converges
\begin{gather*}
R_j\,{\rm d}v_j\to R_X \cH^n+2\theta_x \cH^{n-1}_\Sigma+2\pi(1-r_x)\cH^{n-2}_S  ,
\end{gather*}
where $R_X\colon X\to {\mathbb R}$ is a locally $L^1$ function which is bounded from below, $\Sigma$ is an $n-1$ rectifiable set with $\theta_x\leq 1$, and $2\pi(1-r_x)\cH^{n-2}_S$ is an $n-2$-rectifiable measure supported on the top stratum of the singular set $\operatorname{Sing}(X)$.
\end{Conjecture}

Let us consider two basic examples of the above:
\begin{Example}
Let $X=C\big(S^1_r\big)$ be a metric cone over a circle of radius $r<1$, thus one can picture $X$ as an ice cream cone.  We can smooth~$X$ to nonnegatively curved manifolds~$M_i$ by rounding off the cone tip, see~\cite{LiNa}, and by Gauss--Bonnet one has $R_i\,{\rm d}v_i\to 2\pi(1-r) \delta_0$, where $\delta_0$ is the dirac delta measure at the cone point.
\end{Example}

\begin{Example}
Let $\tilde X=D(0,1)$ be a disk of radius $1$ and let $X$ be its doubling.  Thus $X$ is a~nonsmooth space whose tangent cones are everywhere ${\mathbb R}^2$.   We can smooth $X$ to nonnegatively curved manifolds~$M_i$ by rounding off the boundary, and by Gauss--Bonnet one then has \mbox{$R_i\,{\rm d}v_i\to 2\cH^1_\Sigma$}, where $\Sigma$ is the circle of radius~$1$.
\end{Example}

There would be an interesting corollary of the above Conjecture~\ref{con:2.20}. We had mentioned previously that $\cH^{n-2}(\cS(X))$ can be infinite, however when studying the quantitative stratification $\cS^k_\epsilon$~\cite{ChJiNa} we do have finiteness results on the measure. The following corollary of the above conjecture would give effective understanding of this finiteness for the top stratum (see also~\cite{LiNa} for a more general version of the below in the Alexandrov case):

\begin{Corollary}
	$\cH^{n-2}\big(\cS_\epsilon(X)\cap B_1\big)\leq C(n,v,\lambda) \epsilon^{-1}$. If we denote $\epsilon_i\equiv 2^{-i}$ then we have the Dini type estimate
\begin{gather*}
\sum_i \epsilon_i \cH^{n-2}\big(\big(\cS_{\epsilon_i}\setminus \cS_{\epsilon_{i+1}}\big)\cap B_1\big) \leq C(n,v,\lambda) .
\end{gather*}
\end{Corollary}

\section{Lower Ricci curvature and collapsing}

We now turn our attention to sequences of manifolds with lower Ricci curvature bounds which are collapsing. In order to do this, it becomes particularly important to associate a new piece of information to the limiting process, namely a measure. As the volume of the sequence is tending to zero, one associates with a pointed manifold $\big(M^n_i,g_i,p_i\big)$ the normalized measure $\nu_i\equiv \Vol(B_1(p_i))^{-1}\,{\rm d}v_i$. Then we can consider the measured Gromov--Hausforff limits
\begin{gather*}%\label{e:lowerRicciCollapsing}
\big(M^n_i,g_i,\nu_i,p_i\big)\to (X,d,\nu,p) \qquad \text{s.t.} \quad {\rm Ric}_i\geq -\lambda \quad \text{and} \quad \Vol(B_1(p_i))\to 0 .	
\end{gather*}

It is first worth emphasizing the importance of even being able to study such a sequence. The noncollapsing lower volume bound in the previous section is analogous to an upper bound on the energy in the context of other nonlinear equations, for instance Yang--Mills and nonlinear harmonic maps. In other contexts there are essentially no, or certainly at least very limited, situations in which one can study sequences with energy blowing up, and no general theory. The distinction in the context of lower Ricci curvature is that there is another rigidity which may be exploited in order to produce regularity, namely the splitting theorem.

Currently, the most complete structural theorem about collapsed limits is obtained by combining the results of~\cite{ChC2} and~\cite{CoNa1}:

\begin{Theorem}[Cheeger--Colding \cite{ChC2}, Colding--Naber \cite{CoNa1}]\label{t:structure_collapsing}
Let $\big(M^n_i,g_i,\nu_i,p_i\big)\to (X,d,\nu,p)$ with ${\rm Ric}\geq -\lambda$ and $\Vol(B_1(p_i))\to 0$. Then there exists a unique $k\leq n$ such that $X$ is $k$-rectifiable. In particular, there is a $\nu$-full measure set $\cR^k(X)\subseteq X$ s.t. the tangent cones of $x\in \cR^k$ are unique and isometric to~${\mathbb R}^k$. Further, $\nu\cap \cR^k$ is absolutely continuous with respect to the Hausdorff measure $\cH^k\cap \cR^k$.
\end{Theorem}

The above allows us to associate a unique dimension to a limit space $X$. It is not clear however that this dimension agrees with the Hausdorff dimension. Interestingly, along the regular set $\nu$ and $\cH^k$ are absolutely continuous, and thus the issue is due to the singular set $\cS(X)$, which is a $\nu$-measure zero set:

\begin{oproblem}
	Show the Hausdorff dim of $X$ is same as the rectifiable dimension $k$ from Theorem~{\rm \ref{t:structure_collapsing}}.
\end{oproblem}

The proof of Theorem \ref{t:structure_collapsing} has two distinct parts. In \cite{ChC2} it is first proved that $X=\bigcup_{k} X_k$ decomposes into pieces which are each $k$-rectifiable. In order to prove the uniqueness of $k$, it is then shown in \cite{CoNa1} that tangent cones along geodesics change at a H\"older continuous rate. By then finding geodesics which must intersect different $X_k$ one eventually concludes a contradiction. The H\"older rate is sharp when comparing how the geometries along a geodesic change. Conjecturally, the volume ratio should be behaving even better:

\begin{Conjecture}Let $\gamma\colon (-2,2)\to X$ be a minimizing geodesic, then for $s,t\in (-1,1)$ we have that
\begin{enumerate}\itemsep=0pt
\item[$1)$] $($Lipschitz rate$)$ 	$\frac{\nu(B_r(\gamma(s)))}{\nu(B_r(\gamma(t)))}\leq C(n,\lambda)|t-s|$,
\item[$2)$] $($constant density$)$ if $X$ is noncollapsed, then $\lim\limits_{r\to 0}\frac{\nu(B_r(\gamma(s)))}{\nu(B_r(\gamma(t)))} = 1$.	
\end{enumerate}
\end{Conjecture}

Theorem \ref{t:structure_collapsing} gives a lot in terms of the analytical structure of $X$, but it lacks almost completely a topological understanding of $X$. Maybe the most important open question in this direction is the following:

\begin{oproblem}
Let $\big(M^n_i,g_i,\nu_i,p_i\big)\to (X,d,\nu,p)$ with ${\rm Ric}\geq -\lambda$ and $\Vol(B_1(p_i))\to 0$, then is there an open subset of full $\nu$-measure	$\cR(X)\subseteq X$ which is homeomorphic to a $k$-manifold?
\end{oproblem}

Finally, let us end by discussing the regularity of manifolds with lower Ricci bounds, potentially including those with small volume. The main conjecture out there is one due to Yau, and we state a local version of it below:

\begin{Conjecture}[local Yau conjecture]\label{conj:yau}
	Let $\big(M^n,g\big)$ satisfy ${\rm Ric}\geq -\lambda$, then $\int_{B_1(p)} R \leq C(n,\lambda)$.
\end{Conjecture}

We discussed refinements of the above in the noncollapsing case in Section~\ref{ss:energy_identity}.

\subsection*{Acknowledgements}

The second author was partially supported by the National Science Foundation Grant No.~DMS-1809011.

\pdfbookmark[1]{References}{ref}
\LastPageEnding

\end{document}